\documentclass[a4paper,12pt]{amsart}
\title[A note on Bogomolov-Gieseker type inequality]{{\bf A note on 
Bogomolov-Gieseker type inequality for Calabi-Yau 3-folds}}
\date{}
\author{Yukinobu Toda}

\usepackage{amscd}
\usepackage{amsmath}
\usepackage{amssymb}
\usepackage{amsthm}
\usepackage{float}
\usepackage[dvips]{graphicx}

\usepackage[all,ps,dvips]{xy}

\usepackage{array}
\usepackage{amscd}
\usepackage[all]{xy}

\DeclareFontFamily{U}{rsfs}{%
\skewchar\font127}
\DeclareFontShape{U}{rsfs}{m}{n}{%
<-6>rsfs5<6-8.5>rsfs7<8.5->rsfs10}{}
\DeclareSymbolFont{rsfs}{U}{rsfs}{m}{n}
\DeclareSymbolFontAlphabet
{\mathrsfs}{rsfs}
\DeclareRobustCommand*\rsfs{%
\@fontswitch\relax\mathrsfs}

\theoremstyle{plain}
\newtheorem{thm}{Theorem}[section]

\newtheorem{lem}[thm]{Lemma}

\newtheorem{rmk}[thm]{Remark}
\newtheorem{cor}[thm]{Corollary}

\newtheorem{case}{Case}

\newtheorem{prop-defi}[thm]{Proposition-Definition}
\newtheorem{thm-defi}[thm]{Theorem-Definition}
\newtheorem{lem-defi}[thm]{Lemma-Definition}

\newtheorem{conj}[thm]{Conjecture}

\newdimen\argwidth
\def\db[#1\db]{
 \setbox0=\hbox{$#1$}\argwidth=\wd0
 \setbox0=\hbox{$\left[\box0\right]$}
  \advance\argwidth by -\wd0
 \left[\kern.3\argwidth\box0 \kern.3\argwidth\right]}

\newcommand{\bB}{\mathcal{B}}

\newcommand{\oO}{\mathcal{O}}

\newcommand{\Supp}{\mathop{\rm Supp}\nolimits}
\newcommand{\Hom}{\mathop{\rm Hom}\nolimits}

\newcommand{\Pic}{\mathop{\rm Pic}\nolimits}

\newcommand{\ch}{\mathop{\rm ch}\nolimits}

\newcommand{\Ext}{\mathop{\rm Ext}\nolimits}

\newcommand{\rank}{\mathop{\rm rank}\nolimits}
\newcommand{\Coh}{\mathop{\rm Coh}\nolimits}

\newcommand{\cneq}{\mathrel{\raise.095ex\hbox{:}\mkern-4.2mu=}}
\newcommand{\eqcn}{\mathrel{=\mkern-4.5mu\raise.095ex\hbox{:}}}

\newcommand{\ext}{\mathop{\rm ext}\nolimits}

\newcommand{\length}{\mathop{\rm length}\nolimits}

\begin{document}
\maketitle
\begin{abstract}
The conjectural Bogomolov-Gieseker (BG) type inequality 
for tilt semistable objects on projective
3-folds was proposed 
by Bayer, Macri and the author. In this note, we 
prove our conjecture
for slope stable sheaves
with the smallest first Chern class 
 on certain Calabi-Yau 3-folds, 
e.g. quintic 3-folds.
\end{abstract}
\section{Introduction}
\subsection{Motivation and result}
Let $X$ be a smooth projective 
3-fold over $\mathbb{C}$. 
Given an element 
\begin{align*}
B+i\omega \in H^2(X, \mathbb{C})
\end{align*}
with $\omega$ ample, the heart of a 
bounded t-structure 
\begin{align*}
\bB_{B, \omega} \subset D^b \Coh(X)
\end{align*}
was constructed in~\cite{BMT}, 
following the construction of 
Bridgeland's
stability conditions on projective surfaces~\cite{Brs2}, \cite{AB}. 
The notion of \textit{tilt stability} 
on $\bB_{B, \omega}$ 
was introduced in~\cite{BMT}, and a conjectural 
Bogomolov-Gieseker (BG) type inequality among Chern 
characters of tilt 
semistable objects in $\bB_{B, \omega}$
 was proposed in~\cite[Conjecture~1.3.1]{BMT}. 
Our conjecture in~\cite{BMT} turned out to 
imply several very important results:
construction of Bridgeland stability on projective
3-folds~\cite{BMT}, 
Fujita's conjecture in birational geometry~\cite{BBMT}, 
and Ooguri-Strominger-Vafa's conjecture in string 
theory~\cite{TodBG}. 
In this note, we report a partial progress 
toward the conjectural 
BG type inequality in~\cite{BMT}. 

When $B=0$, the first Chern class on the heart 
$\bB_{0, \omega}$ is always non-negative, and 
plays a role of the rank on the 
category of coherent sheaves. 
A hopeful approach
toward the proof of the main conjecture in~\cite{BMT}
is to use the induction argument on the first
Chern classes of tilt semistable objects, 
as in the proof of BG type inequality
without $\ch_3$. 
(cf.~\cite[Theorem~7.3]{BMT}.) 
As a first step for this induction argument, 
the conjectural BG type inequality 
should be solved when the tilt semistable object has the 
smallest first Chern class. 
In this case, the required result is 
formulated in the following conjecture 
for slope stable sheaves
(cf.~\cite[Conjecture~7.2.3]{BMT}):
\footnote{The statement of~\cite[Conjecture~7.2.3]{BMT}
was more general than Conjecture~\ref{conj:intro} and the 
formulation is slightly different. 
When $\Pic(X)$ is generated by one element, they are obviously 
equivalent.}
\begin{conj}\label{conj:intro}
Let $X$ be a smooth projective 3-fold such that 
$\Pic(X)$ is generated by $\oO_X(H)$ for an ample divisor 
$H$ in $X$. Then for any torsion free slope stable sheaf
$E$ on $X$ with $c_1(E)=[H]$ and $\ch_2(E) H >0$, we have 
the inequality, 
\begin{align}\label{ineq:conj}
\ch_3(E) \le \frac{\ch_2(E) H}{3\ch_0(E)}. 
\end{align}
\end{conj}
The above conjecture was studied 
in~\cite[Example~7.2.4]{BMT}
for rank one torsion free sheaves.
In this case,
the inequality (\ref{ineq:conj}) is 
reduced to Castelnuovo
type inequality
for low degree curves in $X$.  
On the other hand, 
 the higher rank case was not studied in~\cite{BMT}. 
The purpose of this article is to show that,
when $X$ is a  certain Calabi-Yau 3-fold,  
the inequality (\ref{ineq:conj})
is reduced to Castelnuovo type inequality 
even in the higher rank case. 
The main result is as follows: 
\begin{thm}\label{thm:intro}
Let $X$ be a smooth projective Calabi-Yau 3-fold 
such that $\Pic(X)$ is generated by $\oO_X(H)$ for 
an ample divisor $H$ in $X$. 
Suppose that the following inequalities hold:
\begin{align}\label{assum1}
\dim \lvert H \rvert &\ge \frac{7}{6}H^3 -3, \\
\label{assum2}
\chi(\oO_C) &\ge \frac{1}{6}H^3 -C \cdot H,
\end{align}
for any one dimensional subscheme 
$C \subset X$
with $C \cdot H < H^3/2$. Then 
$X$ satisfies Conjecture~\ref{conj:intro}.
Furthermore the inequality (\ref{ineq:conj}) is 
an equality 
only when $E=\oO_X(H)$. 
\end{thm}
As we discussed in~\cite[Example~7.2.4]{BMT}, 
a typical (and important)
example satisfying the assumption of 
Theorem~\ref{thm:intro} is a quintic 3-fold 
in $\mathbb{P}^4$. 
Therefore we obtain the following corollary: 
\begin{cor}\label{cor:intro}
Let $X \subset \mathbb{P}^4$ be a smooth 
quintic 3-fold. Then $X$ satisfies Conjecture~\ref{conj:intro}. 
\end{cor}
In the case of quintic 3-folds, the 
conditions $c_1(E)=[H]$, $\ch_2(E) H >0$
and the Bogomolov-Gieseker inequality~\cite{Bog}, \cite{Gie}
restrict the rank of $E$ up to five.
So a priori, 
the sheaf $E$ could be $\rank(E) \ge 2$.
On the other hand, we do not know any example of
such a sheaf $E$ 
with $\rank(E)\ge 2$. (cf.~Remark~\ref{rmk:possible}.)
The result of Corollary~\ref{cor:intro}
means that $\ch_3(E)$ should obey
 the desired inequality (\ref{ineq:conj}), 
if such a sheaf $E$ exists. 

In general, the
 third Chern character $\ch_3(E)$ is known 
to be bounded by a certain polynomial
 of $\ch_0(E)$, $\ch_1(E)$ and $\ch_2(E)$, 
see~\cite{Langer2}. 
However the evaluation in~\cite{Langer2} is not 
strict to show the inequality (\ref{ineq:conj}).
Although the hypersurface restriction of $E$ plays 
an important role in~\cite{Langer2}, we 
do not take the hypersurface restriction. Instead 
we take the universal extension and 
the classical Bogomolov-Gieseker inequality 
to evaluate the dimensions of cohomology groups.
As far as the author knows, such a method is 
not seen in literatures.

\subsection{Acknowledgement}
This work is supported by World Premier 
International Research Center Initiative
(WPI initiative), MEXT, Japan. This work is also supported by Grant-in Aid
for Scientific Research grant (22684002), 
and partly (S-19104002),
from the Ministry of Education, Culture,
Sports, Science and Technology, Japan.

\subsection{Notation and convention}
In this note, all the varieties are defined over 
$\mathbb{C}$. For $E, F \in \Coh(X)$, 
we denote $h^i(E) \cneq \dim H^i(X, E)$, 
$\mathrm{ext}^i(E, F) \cneq \dim \Ext^i(E, F)$
and $r(E) \cneq \rank(E)$. 
We say $X$ is a Calabi-Yau 3-fold if 
$\dim X=3$, its canonical line bundle is trivial 
and $h^1(\oO_X)=0$. For an ample divisor $H$ 
in a 3-fold $X$ and a torsion free sheaf $E$ on $X$, its 
slope is denoted by 
\begin{align*}
\mu_{H}(E) \cneq \frac{c_1(E) H^2}{r(E)}. 
\end{align*}
The notion of slope stability is defined in the usual way. 
(cf.~\cite{Hu}.)
For a subscheme $Z \subset X$, the defining ideal sheaf 
of $Z$ is denoted by $I_Z$. 

\section{Proof of Theorem~\ref{thm:intro}}
\subsection{Some lemmas}
The key ingredient for the proof of Theorem~\ref{thm:intro}
is the following two lemmas, which may be well-known. 
For the lack of references,  
we give the proofs.  
\begin{lem}\label{lem:1}
Let $X$ be a smooth projective Calabi-Yau 3-fold 
such that $\Pic(X)$ is generated by $\oO_X(H)$ for an 
ample divisor $H$ in $X$. 
Let $E$ be a torsion free slope stable sheaf with 
$c_1(E)=[H]$, and set $V \cneq \Ext^1(E, \oO_X)^{\vee}$. 
Then if we  
take the universal extension
\begin{align}\label{seq:0}
0 \to V \otimes_{\mathbb{C}} \oO_X \to E' \to E \to 0
\end{align}
the sheaf $E'$ is also slope stable. 
\end{lem}
\begin{proof}
We prove the assertion by the induction on $\ext^1(E, \oO_X)$. 
When $\ext^1(E, \oO_X)=0$, then the assertion is obvious. 

Suppose that $\ext^1(E, \oO_X)>0$, and take a non-zero 
element $a \in \Ext^1(E, \oO_X)$. The element $a$
corresponds to the extension, 
\begin{align}\label{ext:a}
0 \to \oO_X \to E_a \to E \to 0. 
\end{align}
We show that $E_a$ is slope stable. 
Suppose by contradiction that $E_a$ is not slope stable. 
Then there is a saturated subsheaf $F \subset E_a$
such that $F$ is slope stable and 
\begin{align}\label{mur}
\mu_H(F) \ge \mu_H(E_a), \quad
r(F)< r(E_a)=r(E)+1.
\end{align}
If we write $c_1(F)=k[H]$, then $k\ge 1$, hence 
$\Hom(F, \oO_X)=0$. It follows that 
the composition 
\begin{align}\label{compose}
F \subset E_a \twoheadrightarrow E
\end{align}
is non-zero, which implies 
$\mu_H(F) \le \mu_H(E)$. Combined with (\ref{mur}),
 we obtain the inequality
\begin{align*}
\frac{H^3}{r(E)+1} \le \frac{kH^3}{r(F)} \le 
\frac{H^3}{r(E)}.
\end{align*}
The above inequality immediately implies $k=1$ and $r(F)=r(E)$. 
Then $\mu_H(F)=\mu_H(E)$, and since $F$ and $E$ are slope 
stable with the same slope, the non-zero morphism (\ref{compose})
is an isomorphism. However this contradicts to that the sequence (\ref{ext:a}) 
is non-split. 

Let $V_a$ be the $\mathbb{C}$-vector space
 $\Ext^1(E_a, \oO_X)^{\vee}$ and take the universal extension, 
\begin{align}\label{seq:2}
0 \to V_a \otimes_{\mathbb{C}} \oO_X \to E_a' \to E_a \to 0. 
\end{align}
Applying $\Hom(-, \oO_X)$ to the sequence (\ref{ext:a}),
we see that 
\begin{align*}
\ext^1(E_a, \oO_X)=\ext^1(E, \oO_X)-1. 
\end{align*}
Hence $E_a'$ is slope 
stable by the assumption of the induction. 
On the other hand, 
composing the sequence (\ref{ext:a}) with (\ref{seq:2}), 
we obtain the exact sequence
\begin{align*}
0 \to V' \otimes_{\mathbb{C}}\oO_X \to E_a' \to E \to 0
\end{align*} 
where $V'$ is a $\mathbb{C}$-vector space 
with $\dim V'=\dim V$. It is easy to see that the above sequence 
is identified with the sequence (\ref{seq:0}), hence 
$E' \cong E_a'$ is slope stable. 
\end{proof}

\begin{lem}\label{lem:2}
In the situation of Lemma~\ref{lem:1}, suppose that 
$r(E)\ge 2$ and 
there is a non-zero element $s\in H^0(X, E)$. Then 
for the associated exact sequence
\begin{align*}
0 \to \oO_X \stackrel{s}{\to} E \to F \to 0
\end{align*}
the sheaf $F$ is also slope stable. 
\end{lem}
\begin{proof}
We first show that $F$ is torsion free. 
If $F$ has a torsion, there is an exact sequence
\begin{align}\label{OAT}
0 \to \oO_X \to A \to T \to 0
\end{align}
where $T$ is a non-zero torsion sheaf
and 
$A \subset E$ is a rank one torsion free sheaf. 
If $\dim \Supp(T)=2$, then $c_1(A)=k[H]$ with $k\ge 1$, 
which contradicts to that $E$ is slope stable. 
Therefore $\dim \Supp(T) \le 1$, hence 
\begin{align*}
\Ext^1(T, \oO_X) &\cong H^2(X, T)^{\vee} \\
&\cong 0. 
\end{align*}
Therefore the sequence (\ref{OAT}) splits,  
which contradicts to that $A$ is torsion free. 

Next suppose that $F$ is not slope stable. Then there 
is a slope stable sheaf $G$ and a surjection 
$F \twoheadrightarrow G$ satisfying 
\begin{align*}
\mu_H(G) \le \mu_H(F), \quad r(G)< r(F)=r(E)-1. 
\end{align*}
Also since there is a surjection 
$E \twoheadrightarrow F \twoheadrightarrow G$
and $E, G$ are slope stable, we have 
$\mu_H(E) \le \mu_H(G)$. 
Hence if we write $c_1(G)=k[H]$, we obtain the
inequality
\begin{align*}
\frac{H^3}{r(E)} \le \frac{kH^3}{r(G)}
\le \frac{H^3}{r(E)-1}. 
\end{align*}
It is immediate to see that there is no solution 
$(k, r(G))$ satisfying the above inequality and $r(G)< r(E)-1$. 
Hence $F$ is slope stable. 
\end{proof}
As a corollary of Lemma~\ref{lem:2}, we have the following: 
\begin{cor}\label{cor}
In the situation of Lemma~\ref{lem:1}, there is 
an exact sequence of the form
\begin{align}\label{seq:3.5}
0 \to \oO_X^{\oplus m} \to E \to F \to 0
\end{align}
such that $F$ is either a rank one torsion free sheaf 
or a slope stable sheaf with $r(F) \ge 2$ and $h^0(F)=0$. 
\end{cor}
\begin{proof}
We show the assertion by the induction of $\theta(E)$ defined by 
\begin{align*}
\theta(E) \cneq \mathrm{min}\left\{ h^0(E), r(E)-1 \right\}. 
\end{align*}
The assertion is obvious when $\theta(E)=0$. 
Suppose that $\theta(E)>0$, i.e. $h^0(E) \neq 0$ and 
$r(E)\ge 2$.  Then 
there is a non-zero element $s\in H^0(X, E)$. 
If we take the exact sequence
\begin{align}\label{seq:3}
0 \to \oO_X \stackrel{s}{\to} E \to F_s \to 0 
\end{align}
then $F_s$ is slope stable by Lemma~\ref{lem:2}. 
By applying 
$\Hom(\oO_X, -)$ to the sequence (\ref{seq:3}), 
we see $h^0(F_s)=h^0(E)-1$. 
Hence we have 
$\theta(F_s)=\theta(E)-1$, 
and by the assumption of the induction,  
 there is an exact
sequence
\begin{align}\label{seq:4}
0 \to \oO_X^{\oplus m'} \to F_s \to F \to 0
\end{align}
such that $F$ is a rank one torsion free sheaf or a
slope stable sheaf with 
$r(F) \ge 2$ and $h^0(F)=0$. 
The desired exact sequence (\ref{seq:3.5}) is obtained by 
combining the sequence (\ref{seq:4}) with (\ref{seq:3}). 
\end{proof}

\subsection{Proof of Theorem~\ref{thm:intro}}
\begin{proof}
Let $X$
 be as in the statement of Theorem~\ref{thm:intro}, 
and $E$ a slope stable sheaf on $X$ with $c_1(E)=[H]$
and $\ch_2(E) H >0$. 
By Corollary~\ref{cor}, there is an exact sequence of the form
\begin{align}\label{seq:OEF}
0 \to \oO_X^{\oplus m} \to E \to F \to 0
\end{align}
such that either $F$ is a rank one torsion free sheaf or 
a slope stable sheaf with $r(F) \ge 2$ and $h^0(F)=0$. 
Note that in the first case, we have 
$F \cong \oO_X(H) \otimes I_Z$
for a subscheme $Z \subset X$ with $\dim Z \le 1$. 
We evaluate $\ch_3(E)=\ch_3(F)$ 
by dividing into the following three cases: 
\begin{case}
We have 
$F \cong \oO_X(H) \otimes I_Z$ with $\dim Z=0$. 
\end{case}
In this case, we have 
\begin{align*}
\Ext^1(F, \oO_X) \cong H^2(X, F)^{\vee}
\end{align*}
by the Serre duality, which is zero by 
the cohomology exact sequence associated to 
the sequence
\begin{align*}
0 \to \oO_X(H) \otimes I_Z \to \oO_X(H) \to \oO_Z \to 0, 
\end{align*}
and the Kodaira vanishing $h^2(\oO_X(H))=0$. 
Hence the sequence (\ref{seq:OEF}) splits if $m>0$, 
which contradicts to the slope stability of $E$. 
Therefore $E \cong \oO_X(H) \otimes I_Z$, 
and 
\begin{align*}
\ch_3(E)&=\frac{1}{6}H^3 - \length \oO_Z, \\
\frac{\ch_2(E)}{3\ch_0(E)}&= \frac{1}{6}H^3 - \frac{1}{3}
\length \oO_Z.
\end{align*}
The above equalities imply 
 the inequality (\ref{ineq:conj}), and the 
equality holds only when $Z=\emptyset$. 

\begin{case}
We have 
$F \cong \oO_X(H) \otimes I_Z$ with $\dim Z=1$. 
\end{case}
In this case, 
$\ch_2(E) H =\ch_2(F) H>0$ is equivalent to 
\begin{align*}
0< H \cdot Z < \frac{1}{2}H^3. 
\end{align*}
Applying 
the assumption (\ref{assum2}), we have 
\begin{align*}
\ch_3(E) &=\ch_3(F) \\
&= \frac{1}{6}H^3 - H \cdot [Z] -\chi(\oO_Z) \\
&\le 0. 
\end{align*}
Therefore the inequality (\ref{ineq:conj}) holds. 

\begin{case}
We have $r(F)\ge 2$ and $h^0(F)=0$. 
\end{case}
By the Riemann-Roch theorem, we have 
\begin{align}\notag
\chi(E) &\cneq \sum_{i\ge 0}
(-1)^i h^i(E) \\
\notag
&= \ch_3(E) + \frac{1}{12}H \cdot c_2(X) \\
\label{c3:eq}
&= \ch_3(E) + \dim \lvert H \rvert -\frac{1}{6}H^3 +1. 
\end{align}
Here the last equality follows from the Riemann-Roch theorem 
applied for $\oO_X(H)$ and the Kodaira vanishing $h^i(\oO_X(H))=0$
for $i>0$. 
On the other hand, since $\chi(\oO_X)=0$ and $h^0(F)=0$, we have 
\begin{align}\label{c3:ineq}
\chi(E)= \chi(F) \le h^2(F)=\ext^1(F, \oO_X).
\end{align}
Here the last equality follows from the Serre duality. 
Let us take the universal extension
\begin{align*}
0\to \oO_X^{\oplus \ext^1(F, \oO_X)} \to F' \to F \to 0. 
\end{align*} 
By Lemma~\ref{lem:1}, the sheaf $F'$ is slope stable 
satisfying
\begin{align*}
\ch_0(F')&=\ch_0(F)+ \ext^1(F, \oO_X), \\
\ch_i(F')&=\ch_i(F)=\ch_i(E), \ i\ge 1.
\end{align*}
By the Bogomolov-Gieseker inequality~\cite{Bog}, \cite{Gie}
applied for $F'$, we obtain the inequality
\begin{align*}
\left(H^2 -2(\ch_0(F)+ \ext^1(F, \oO_X)) \ch_2(E)\right) H \ge 0.
\end{align*}
Since $\ch_2(E) H>0$ by the assumption, the above inequality is 
equivalent to 
\begin{align}\label{c3:ineq2}
\ext^1(F, \oO_X) \le \frac{H^3}{2\ch_2(E) H} -\ch_0(F). 
\end{align}
By (\ref{c3:eq}), (\ref{c3:ineq}), (\ref{c3:ineq2})
and noting $\ch_2(E) H\ge 1/2$, $\ch_0(F)\ge 2$, 
we obtain the inequality
\begin{align}\notag
\ch_3(E) &\le \frac{H^3}{2\ch_2(E) H} -\ch_0(F) + \frac{1}{6}H^3 -\dim \lvert H \rvert -1 \\
\label{strengh}
&\le \frac{7}{6}H^3 -\dim \lvert H \rvert  -3 \\
\notag
&\le 0. 
\end{align}
The last inequality follows from the assumption (\ref{assum1}).
Therefore the inequality (\ref{ineq:conj}) holds.  
\end{proof}

\begin{rmk}
When $H^3$ is even, we
have $\ch_2(E) H \ge 1$, hence we 
have the better inequality in (\ref{strengh}). 
Consequently, 
we can weaken
the assumption (\ref{assum1}) to be 
\begin{align*}
\dim \lvert H \rvert \ge \frac{2}{3}H^3 -3. 
\end{align*}
For instance, the complete intersections 
\begin{align*}
(2, 4) \subset \mathbb{P}^5, \ (2, 2, 3) \subset \mathbb{P}^6,
\end{align*} 
satisfy the above inequality. However the Castelnuovo type 
inequality (\ref{assum2}) remains open, 
as we discussed in~\cite[Example~7.2.4]{BMT}.
\end{rmk}

\begin{rmk}\label{rmk:possible}
Let $X$ be a quintic 3-fold
and $E$ a slope stable sheaf $E$ on $X$ satisfying 
the assumption of Conjecture~\ref{conj:intro}. 
Then 
the Bogomolov-Gieseker inequality restricts the 
possibility for the vector $(r(E), c_2(E) H)$
to be one of the following: 
\begin{align*}
(1, 0), \ (1, 1), \ (1, 2), \ (2, 2), \ (3, 2), \ (4, 2), \ (5, 2). 
\end{align*}
However we don't know any example of $E$ with 
$r(E)\ge 2$ and $c_2(E) H =2$. 
\end{rmk}

\begin{rmk}
The result of Corollary~\ref{cor:intro} easily 
shows that the conjecture in~\cite[Conjecture~1.3.1]{BMT}
is true if $X$ is a quintic 3-fold, $B=0$ and $c_1(E)=[H]$. 
\end{rmk}

\begin{rmk}\label{rmk:ch3}
When $X$ is a quintic 3-fold, the proof of Theorem~\ref{thm:intro}
shows that $\ch_3(E) \le 0$ except $E=\oO_X(H)$. 
\end{rmk}

\begin{rmk}
The exact sequences of sheaves (\ref{seq:0}), (\ref{seq:3.5})
are nothing but exact sequences
\begin{align*}
&0 \to E' \to E \to V\otimes_{\mathbb{C}} \oO_X[1] \to 0, \\
&0\to E \to F \to \oO_X^{\oplus m}[1] \to 0, 
\end{align*}
in the tilted heart 
$\bB_{0, H} \subset D^b \Coh(X)$. 
This observation may be helpful for 
the study of~\cite[Conjecture~1.3.1]{BMT}
when the first Chern class of $E$ 
is higher than $[H]$. 
\end{rmk}
\begin{rmk}\label{rmk:ind}
The result of Corollary~\ref{cor:intro} 
also implies a partial result on~\cite[Conjecture~1.3.1]{BMT}
when $X$ is a quintic 3-fold, $B=0$ and $c_1(E)=2[H]$. 
Indeed, in the notation of~\cite{BMT},
suppose that $E \in \bB_{0, H}$ satisfies $c_1(E)=2[H]$
and $\nu_{t[H]}(E)=0$ for some 
$t \in \mathbb{R}_{>0}$. If  
$E$ is not (essentially) a 
sheaf, then $E$ fits into an exact sequence in $\bB_{0, H}$
\begin{align*}
0 \to F_1[1] \to E \to F_2 \to 0
\end{align*}
such that $F_i$ are 
slope stable sheaves with $c_1(F_1)=[H]$, $c_1(F_2)=-[H]$.
The $\nu_{t[H]}$-stability 
implies 
\begin{align*}
\nu_{t[H]}(F_1[1]) \le \nu_{t[H]}(E)=0 \le \nu_{t[H]}(F_2),
\end{align*}
which shows $\ch_2(F_i) H>0$. 
Hence we can apply Corollary~\ref{cor:intro} for each $F_i$
(and noting Remark~\ref{rmk:ch3})
to conclude that $\ch_3(E) \le 5/6$. 
This is enough to show 
the inequality conjectured
 in~\cite[Conjecture~1.3.1]{BMT}.
We hope that the above argument is 
generalized to the case that 
$c_1(E)$ is higher than $2[H]$ 
as an induction argument.
\end{rmk}
\begin{rmk}
The argument of Remark~\ref{rmk:ind} shows that, 
 in the case of $c_1(E)=2[H]$, the conjecture in~\cite[Conjecture~1.3.1]{BMT} is reduced to the case that $E$ is a sheaf, 
e.g. $E$ is a rank one torsion free sheaf. 
However in this case,
we need a 
Castelnuovo type inequality, stronger than the known one,  
to show the desired inequality in~\cite[Conjecture~1.3.1]{BMT}.
\end{rmk}

Institute for the Physics and 
Mathematics of the Universe, 

Todai Institute for Advanced Studies (TODIAS), 
University of Tokyo,

5-1-5 Kashiwanoha, Kashiwa, 277-8583, Japan.

\textit{E-mail address}: yukinobu.toda@ipmu.jp

\end{document}